\documentclass[11pt,leqno]{article}

 \usepackage{geometry}
\usepackage{amssymb}
\usepackage{amsthm,amscd,amsmath}
\usepackage{tikz-cd}
\usepackage{color}
\usetikzlibrary{decorations.pathmorphing}

\numberwithin{equation}{section}
\usepackage{enumerate}
\usepackage{cite}
\usepackage{mathrsfs}
\usepackage[pagebackref]{hyperref}

\newcommand{\bpm}{\begin{pmatrix}}
\newcommand{\epm}{\end{pmatrix}}

\newcommand{\N}{\mathbf{N}}

\newcommand{\cat}[1]{\mathbf{#1}}
\newcommand{\twocat}[1]{\underline{\mathbf{#1}}}

\newcommand{\mc}{\mathcal}

\newcommand{\iso}{\cong}

\newtheorem{thm}[equation]{Theorem}
\newtheorem{cor}[equation]{Corollary}
\newtheorem{lemma}[equation]{Lemma}

\theoremstyle{definition}
\newtheorem{definition}[equation]{Definition}

\theoremstyle{remark}

\title{On Whitehead's theorem beyond pointed connected spaces}
\author{Kevin Arlin\footnote{The author's last name was previously Carlson.}}

\begin{document}
\maketitle

The aim of this note is to prove a Whitehead theorem for spaces, specifically, 
Theorem \ref{whiteheadforspaces}. 
Of course, the classical Whitehead theorem, in the form that
 a map $f:X\to Y$ of spaces is a  homotopy equivalence if it induces 
 isomorphisms on all homotopy groups at all base points, is already
 about spaces. We are going, instead, for the stronger form the Whitehead 
 theorem takes 
in pointed connected spaces. Let $\cat{Hot}$ denote the category 
of CW complexes (equivalently, Kan complexes) and homotopy classes of continuous 
maps. Let $\cat{Hot}_{*,c}$ denote the category of pointed connected 
CW complexes and equivalence classes of pointed maps up to homotopies through
pointed maps. Then the Whitehead theorem \cite{whitehead} 
can be interpreted as follows:

\begin{thm}[Whitehead]\label{WhiteheadClassical}
The category $\cat{Hot}_{*,c}$ admits the set of spheres 
$\{S^n\}$ as a strong generator. 
\end{thm}
That is, the spheres detect isomorphisms in $\cat{Hot}_{*,c}$. 

This version of Whitehead's theorem is at the heart of various aspects of
modern homotopy theory. For instance, the triangulated categories underlying
presentable stable model categories or $\infty$-categories are known to 
admit strong generators. In that context and in spaces, the existence
of a strong generator is crucial to the proof of the Brown representability
theorem. Thus it is an important flaw of the category $\cat{Hot}$ 
that it does \emph{not} admit
a generator. This was claimed in \cite{heller3}, but Heller's proof only
goes through to show that there exists no generator consisting of 
\emph{finite} spaces. 
A correct proof that there exists no generator whatsoever 
will appear in future work of the author and Dan Christensen. 

We will show that, while $\cat{Hot}$ lacks any generator, the 2-category
of spaces $\twocat{Hot}$ admits a very manageable one, namely, the
finite-dimensional tori. By $\twocat{Hot}$, we mean the 2-category
whose objects are CW complexes (again, equivalently, Kan complexes) 
and whose hom-categories are the fundamental
groupoids of mapping spaces, that is, $\twocat{Hot}(X,Y)=\Pi_1(Y^X)$. 
Thus in $\twocat{Hot}(X,Y)$, an object is a map $f:X\to Y$ while a morphism
is an equivalence class of homotopies $X\times I\to Y$ between $f$ and $g$
up to homotopies $X\times I\times I\to Y$ which are constant at $f$ and $g$
on $X\times \{0\}\times I$, respectively, $X\times \{1\}\times I$. 

By a generator for a 2-category, we mean the following:
\begin{definition}\label{2-strong generator}
A small set $\mc G$ of objects in a 2-category $\mc K$ 	
constitutes a strong generator for $\mc K$ if $\mc G$ jointly
detects equivalences. That is, whenever $f:X\to Y$ is a morphism in
$\mc K$ such that, for every $G\in \mc G$, 
$\mc K(G,f):\mc K(G,X)\to \mc K(G,Y)$ is an equivalence of categories,
we may conclude that $f$ is an equivalence in $\mc K$. 
\end{definition}

We now develop the necessary lemmata for Theorem \ref{whiteheadforspaces}. 
First, we must compute the homotopy groups of free loop spaces based at constant
loops. 
For any space $X$, let $LX$ be its free loop space, together with the embedding
$c:X\to LX$ of constant loops. Note that $c$ splits any evaluation map 
$e:LX\to X$. We abuse notation by writing $\pi_n(LX,x)$ for $\pi_n(LX,c(x))$. 
Then for any $x\in X$, since 
$\pi_n(e,x)\circ\pi_n(c,x)=\mathrm{id}_{\pi_n(X,x)}$,  we have a 
decomposition
\begin{equation}\label{splitting}\pi_n(LX,x)\cong \pi_n(X,x)\ltimes\ker(\pi_n(e,x)).\end{equation} 
Now it remains only to
compute the normal factor $\ker(\pi_n(e,x))$. 

For this, we identify the splitting of Equation \ref{splitting} 
with that induced by the long exact sequence of the fibration $e:LX\to X$.
Recall that $e$ has fiber
over $x$ the space $\Omega(X,x)$ of loops based at $x$ and denote by 
$\iota:\Omega(X,x)\to LX$ the inclusion of the fiber. 
Then we have $\ker(\pi_n(e,x))=\mathrm{im}(\pi_n(\iota,x))$. Since 
$e$ is split, the boundary morphisms in the long exact sequence vanish, 
so we have a canonical isomorphism $\ker(\pi_n(e,x))\iso \pi_n(\Omega(X,x))$. 
Composing with the 
usual identification of homotopy groups of based loop spaces, we have proved

\begin{lemma}\label{pinL}
We have isomorphisms $(\pi_n(c,x),\pi_n(\iota,x))$, natural on the homotopy 
category $\cat{Hot}_*$ of pointed spaces, identifying $\pi_n(LX,x)$
as a semidirect product $\pi_n(X,x)\ltimes \pi_{n+1}(X,x)$.
\end{lemma}

Of course, when $n>1$, the semidirect product is actually direct. 
A straightforward induction shows, furthermore, that higher homotopy groups of 
iterated free loop spaces are given by the binomial theorem. 

\begin{cor}\label{pinLn}
We have natural isomorphisms
$\pi_k(L^n X,x)\cong\bigoplus_{i=0}^n \pi_{k+i}(X,x)^{\oplus \binom n i}$
for each $k>1$.
\end{cor}
Finally, we expose a consequence of the naturality of these
isomorphisms which we will need below.
\begin{cor}\label{naturality}
Let $f:X\to Y$ be any map of spaces, $x\in X$, and fix $n\in\N$ with $n\geq 2$. 
Suppose $f$ induces
isomorphisms $\pi_k(f,x):\pi_k(X,x)\to \pi_k(Y,f(x))$ as well as
$\pi_1(L^kf,x):\pi_1(L^{k}X,x)\to \pi_1(L^{k}Y,f(x))$, for each $k<n$. 
Then $\pi_n(f,x):\pi_n(X,x)\to \pi_n(Y,f(x))$ is an isomorphism.
\end{cor}
\begin{proof}
We first apply Lemma \ref{pinL} to obtain the map of short exact
sequences of groups
\[\begin{tikzcd}
				1\ar[r]&\pi_2(L^{n-2}X,x)\ar[d,"p"]\ar[r]&\pi_1(L^{n-1}X,x)\ar[d,"q"]\ar[r]&\pi_1(L^{n-2}X,x)\ar[d,"r"]\ar[r]&1\\
				1\ar[r]&\pi_2(L^{n-2}Y,f(x))\ar[r]&\pi_1(L^{n-1}Y,f(x))\ar[r]&\pi_1(L^{n-2}Y,f(x))\ar[r]&1
\end{tikzcd}\]
in which the vertical arrows are induced by $f$. By assumption, $q$ and 
$r$ are isomorphisms, which implies $p$ is. Define the functor 
$T(X,x)$ on $\cat{Hot}_*$ 
by $T(X,x)=\oplus_{i=0}^{n-3} \pi_{2+i}(X,x)^{\oplus \binom{n-2} i}$.
Corollary \ref{pinLn} identifies
$T$ with a split subfunctor of $\pi_2(L^{n-1}X,x)$,
yielding another morphism of short exact sequences
\[\begin{tikzcd}
				0\ar[r]&\pi_n(X,x)\ar[d,"s"]\ar[r]&\pi_2(L^{n-2}X,c(x))\ar[d,"p"]\ar[r]&T(X,x)\ar[d,"T(f)"]\ar[r]&0\\
				0\ar[r]&\pi_n(Y,f(x))\ar[r]&\pi_2(L^{n-2}Y,c(f(x)))\ar[r]&T(Y,f(x))\ar[r]&0
\end{tikzcd}\]
We have just shown that $p$ is an isomorphism, and by assumption, so is
$T(f)$. Thus $s$ is an isomorphism, as was to be shown. 
\end{proof}

We turn to the main theorem.

\begin{thm}\label{whiteheadforspaces}
Let $T^n$ denote the $n$-torus $\prod_{i=1}^nS^1$. Then the family
 $\{T^n\}_{n\in\N}$ constitutes a strong generator for
the 2-category $\twocat{Hot}$ of spaces. 
\end{thm}
\begin{proof}
Assume that $f:X\to Y$ is such that 
$\twocat{Hot}(T^n,f)$ is an equivalence of groupoids for every
$n\in\N$. We shall show that $f$ induces an isomorphism on connected
components and homotopy groups, 
$\pi_n(f,x):\pi_n(X,x)\cong \pi_n(Y,f(x))$
for every $n$ and for every $x\in X$. 
The $\pi_0$ and $\pi_1$ isomorphisms are immediate from the 
assumption, since the 0-torus $T^0$ is a point and 
$\twocat{Hot}(*,X)=\Pi_1(X)$ is the fundamental groupoid of $X$, while we
have natural isomorphisms $\pi_0\Pi_1(X)\iso \pi_0 X$ and 
$\Pi_1(X)(x,x)\iso \pi_1(X,x)$. 

Now, given $n>1$, we can assume that $f$ induces isomorphisms 
$\pi_k(X,x)\iso \pi_k(Y,f(x))$ for every $k<n$. By assumption, 
$f$ induces an isomorphism between $\twocat{Hot}(T^k,X)(x,x)$ and
$\twocat{Hot}(T^k,Y)(f(x),f(x))$, for every $k$.
That is, we have an isomorphism 
\[\pi_1(L^kf,x):\pi_1(L^kX,x)\iso \pi_1(L^kY,f(x))\] 
Thus the assumptions of Corollary
\ref{naturality} hold, and $\pi_n(f,x):\pi_n(X,x)\to \pi_n(Y,f(x))$ is an 
isomorphism, as was to be shown. 
  \end{proof}

\bibliographystyle{abbrv}
\bibliography{bib.bib}

\end{document}